\documentclass[11pt,a4paper,reqno]{amsart}
\usepackage{amsmath}
\usepackage{amssymb,latexsym}
\usepackage[dvips]{graphicx}
\usepackage{color}

\begin{document}
\medskip
  \title{ Remark on the paper "On the eigenvalues of the Witten-Laplacian on compact Riemannian manifolds " by Q.-M. Cheng and L. Zeng}
\author{  S. Ilias}
\date{13/04/13}

\maketitle
\indent 
 On April 11, 2013, a preprint was posted on arXiv by Q.-M. Cheng and L. Zeng \cite{CZ}, about eigenvalues of the Witten-Laplacian acting on functions. However, it is very well known that the Witten Laplacian $$\Delta_{f}:=\Delta+<\nabla f,.>$$ on a compact measure metric space $(M,g,e^{-f}\,dv )$ is unitarily equivalent to the Schr\"{o}dinger operator $\Delta+\frac{1}{2}\Delta f+\frac{1}{4}|\nabla f|^{2}$ on $(M,g)$ and it has thus the same spectrum (see for instance \cite{S}). Therefore, all the results of this above mentioned preprint are immediate consequences of that concerning the Schr\"{o}dinger operator  $\tilde{\Delta}_{f}:=\Delta+\frac{1}{2}\Delta f+\frac{1}{4}|\nabla f|^{2}$ on $(M,g)$ already obtained in the papers \cite{EI, EHI, IM1, IM2, WX},  provided replacing the potential $q$ or $V$ in that papers by $\frac{1}{2}\Delta f+\frac{1}{4}|\nabla f|^{2}$   (obviousely, this is in particular, the case for the self shrinkers).\\

\indent Using the same observation of unitary equivalence above, we also note that a lot of published results concerning eigenvalues of the Laplace-Witten operator, can be deduced from that concerning Schr\"{o}dinger operators.
 
\textsc{Said Ilias}

{\em Universit\'e F. Rabelais, D\'ep. de Math\'ematiques, Tours

ilias@univ-tours.fr}

\end{document}